\documentclass[preprint,12pt]{elsarticle}

\usepackage{epsfig}

\usepackage{amssymb}
\usepackage{amsmath}

\usepackage{amsthm}
\usepackage{graphicx}

\newtheorem{thm}{Theorem}

\newdefinition{rmk}{Remark}
\newproof{pf}{Proof}
\newproof{pot}{Proof of Theorem \ref{main}}

\begin{document}

\begin{frontmatter}


\ead{aakirac@pau.edu.tr}

\title{On the Ambarzumyan's theorem for the Quasi-periodic Problem}


\author{Alp Arslan K\i ra\c{c}}

\address{Department of Mathematics, Faculty of Arts and Sciences, Pamukkale
University, 20070, Denizli, Turkey}

\begin{abstract}
We obtain the classical Ambarzumyan's theorem for the Sturm-Liouville operators $L_{t}(q)$ with $q\in L^{1}[0,1]$ and quasi-periodic boundary conditions, $t\in [0,2\pi)$, when there is not any additional condition on the potential $q$.
\end{abstract}

\begin{keyword}
Ambarzumyan theorem; inverse spectral theory; Hill operator



\end{keyword}

\end{frontmatter}


\section{Introduction}
In this study we consider the Sturm-Liouville operator $L_{t}(q)$ generated in the space $L^{2}[0, 1]$ by the
expression
\begin{equation}  \label{1}
-y^{\prime\prime}+q(x)y
\end{equation}
and the quasi-periodic boundary conditions
\begin{equation}\label{q-per}
y(1)=e^{it}y(0),\qquad y^{\prime}(1)=e^{it}y^{\prime}(0),
\end{equation}
where $q\in L^{1}[0,1]$ is a real-valued function and $t$ is a fixed real number in $[0,2\pi)$. Note that the operator $L_{t}(q)$ is  self-adjoint and the cases $t=0$ and $t=\pi$ correspond to the periodic and antiperiodic problems, respectively. Since the spectrum $S(L(q))$ of Hill operator $L(q)$ generated in the space $L^{2}(-\infty, \infty)$ by expression (\ref{1}) with periodic potential $q$ is the union of the spectra $S(L_{t} (q))$  of the operators $L_{t}(q)$ for $t\in [0,2\pi)$ (e.g., see \cite{Eastham}), the operators $L_{t}(q)$ have a fundamental role in the spectral
theory of the operator $L(q)$. In $1929$, Ambarzumyan \cite{Ambarz} obtained the following theorem considered as the first theorem in inverse spectral theory:

\emph{If $\{n^{2}: n=0,1,\ldots\}$ is the spectrum of the Sturm-Liouville operator with Neumann boundary condition, then $q=0$ a.e.}

 In \cite{Chern}, Chern and Shen proved Ambarzumyan's theorem for the Sturm-Liouville
differential systems with Neumann boundary conditions. Later, in \cite{corri}, by imposing an additional condition on the potential they extended the classical Ambarzumyan's theorem for the Sturm-Liouville equation to the general separated boundary conditions. See basics and further references in \cite{Levitan, Lieberman}.

At this point we refer in particular to \cite{Ambarzcoupled,anoteinver}. In \cite{Ambarzcoupled}, for the vectorial Sturm-Liouville
problem under periodic or antiperiodic boundary conditions, Yang-Huang-Yang found two analogs of Ambarzumyan's theorem. Their result supplements the P{\"{o}}schel-Trubowitz inverse spectral theory \cite{Poschel}. More recently, Cheng-Wang-Wu \cite {anoteinver} proved the following theorem:

\emph{(a) If all eigenvalues of the operator $L_{0}(q)$ are nonnegative and they include $\{(2n\pi)^{2}: n\in \mathbb{N}\}$, then $q=0$ a.e.}

\emph{(b) If all eigenvalues of the operator $L_{\pi}(q)$  are not less than $\pi^{2}$ and they include $\{(2n\pi-\pi)^{2}: n\in \mathbb{N}\}$, and
\begin{equation}  \label{cond}
\int_{0}^{1}q(x)\,cos(2\pi x)\,dx\geq 0,
\end{equation} then $q=0$ a.e.}

The present work was stimulated by the papers \cite{corri, anoteinver}. For the first time, we obtain Ambarzumyan's theorem for the operator $L_{t}(q)$ with $t\in[0,2\pi)$, generated by quasi-periodic boundary
conditions (\ref{q-per}). The result established below show that the potential $q$ can be determined from one spectrum and there is not any additional condition on $q$ such as (\ref{cond}) for the operator $L_{t}(q)$ with $t=\pi$ (see also \cite{corri, Ambarzcoupled}). The result of this paper is the following.

\begin{thm} \label{main}
If first eigenvalue of the operator $L_{t}(q)$ for any fixed number $t$ in $[0,2\pi)$ is  not less than the value of $\min\{t^{2}, (2\pi-t)^{2}\}$ and the spectrum $S(L_{t} (q))$ contains the set
 $\{(2n\pi -t)^{2}: n\in \mathbb{N}\}$, then $q=0$ a.e.
\end{thm}

\section{Preliminaries and Proof of the result}
 We now introduce some preliminary facts. In \cite{Melda.O} (see also \cite{Veliev:Kira�}), without using the assumption $q_{0}=0$, they proved the following result:

\emph{The eigenvalues $\lambda_{n}(t)$ of the operator $L_{t}(q)$ for $q\in L^{1}[0,1]$ and $t\neq 0,\pi$, satisfy the following asymptotic formula
\begin{equation}  \label{melvel}
\lambda_{n}(t)=(2\pi n+t)^{2}+ q_{0}+O\left(n^{-1}ln|n|\right)\qquad \textrm{as $|n|\rightarrow \infty$},
\end{equation}
where $q_{n}=(q, e^{i2\pi nx})$ for $n\in \mathbb{Z}$ and $(.\, , .)$ is the inner product in $L^{2}[0, 1]$.
}

Note that when $q=0$, $(2\pi n+t)^{2}$ for $n\in \mathbb{Z}$ is the eigenvalue of the operator $L_{t}(0)$ for any fixed $t\in [0,2\pi)$ corresponding to the eigenfunction $e^{i(2\pi n+t)x}$.
\begin{pot}
Using the assumption that, for any $n\in \mathbb{N}$, $(2n\pi -t)^{2}$  belongs to the spectrum $S(L_{t} (q))$ and taking into account that, for sufficiently large $|n|$, the asymptotic formulas (\ref{melvel}) for $t\neq 0,\pi,$ and, in \cite{anoteinver}, (1.2)-(1.3) for $t=0,\pi$ (see Theorem 1.1. of \cite{anoteinver}), we obtain
 \begin{equation}  \label{q0}
q_{0}=\int_{0}^{1}q(x)\,dx= 0.
\end{equation}
Let us show that, for fixed $t\in [0,2\pi)$, the first eigenvalue of the operator $L_{t}(q)$  is either $t^{2}$ or $(2\pi-t)^{2}$ corresponding to the eigenfunctions $y=e^{itx}$ or $y=e^{i(-2\pi+t)x}$, respectively. First, suppose that the value of $\min\{t^{2}, (2\pi-t)^{2}\}$ is $t^{2}$. By the variational principle and (\ref{q0}), we have for $y=e^{itx}$
\begin{equation}  \label{ineq}
t^{2}\leq\lambda_{0}(t)\leq \frac{\int_{0}^{1}-\bar{y}y^{\prime\prime}dx+\int_{0}^{1}q(x)|y|^{2}dx}{(y,y)}=t^{2}+q_{0}=t^{2}.
\end{equation}
This implies that the first eigenvalue of the operator $L_{t}(q)$ is $\lambda_{0}(t)=t^{2}$ and the test function $y=e^{itx}$ is the first eigenfunction of the operator. Thus, Substituting the expressions $y=e^{itx}$ and $\lambda_{0}(t)=t^{2}$ into the equation $$-y^{\prime\prime}+q(x)y=\lambda y,$$
we get $q=0$ in $L^{1}[0,1]$. Similarly, one can readily show that if the value of $\min\{t^{2}, (2\pi-t)^{2}\}$ is $(2\pi-t)^{2}$, then the function $y=e^{i(-2\pi+t)x}$ is the first eigenfunction corresponding to the first eigenvalue $(2\pi-t)^{2}$ and $q=0$ in $L^{1}[0,1]$.\qed
 \end{pot}
 \begin{rmk}
 Note that instead of the subset $\{(2n\pi -t)^{2}: n\in \mathbb{N}\}$ of the spectrum $S(L_{t} (q))$ in Theorem \ref{main} if we use either of the subsets  $$\{(2n\pi +t)^{2}: n\in \mathbb{N}\}\,\, ,\{m^{2}: m\,\textrm{is either $(2n\pi -t)$ or $(2n\pi +t)$ for all $n\in \mathbb{N}$}\},$$ then the assertion of Theorem \ref{main} remains valid.
\end{rmk}


%



\begin{thebibliography}{10}
\expandafter\ifx\csname url\endcsname\relax
  \def\url#1{\texttt{#1}}\fi
\expandafter\ifx\csname urlprefix\endcsname\relax\def\urlprefix{URL }\fi
\expandafter\ifx\csname href\endcsname\relax
  \def\href#1#2{#2} \def\path#1{#1}\fi

\bibitem{Eastham}
M.~S.~P. Eastham, The Spectral Theory of Periodic Differential Operators,
  Scottish Academic Press, Edinburgh, 1973.

\bibitem{Ambarz}
V.~Ambarzumian, {\"{U}}ber eine \textrm{Frage} der \textrm{Eigenwerttheorie},
  Zeitschrift f{\"{u}}r Physik 53 (1929) 690--695.

\bibitem{Chern}
H.~H. Chern, C.~L. Shen, On the n-dimensional \textrm{Ambarzumyan's} theorem,
  Inverse Problems 13~(1) (1997) 15--18.

\bibitem{corri}
H.~H. Chern, C.~K. Lawb, H.~J. Wang, Corrigendum to \textrm{“Extension of
  Ambarzumyan's theorem} to general boundary conditions”, J. Math. Anal. Appl.
  309 (2005) 764--768.

\bibitem{Levitan}
B.~M. Levitan, M.~G. Gasymov, Determination of a differential equation by two
  of its spectra, Usp. Mat. Nauk 19 (1964) 3--63.

\bibitem{Lieberman}
H.~Hochstadt, B.~Lieberman, An inverse sturm-liouville problem with mixed given
  data, SIAM J. Appl. Math. 34 (1978) 676--680.

\bibitem{Ambarzcoupled}
C.~F. Yang, Z.~Y. Huang, X.~P. Yang, Ambarzumyan's theorems for vectorial
  sturm-liouville systems with coupled boundary conditions., Taiwanese J. Math.
  14~(4) (2010) 1429--1437.

\bibitem{anoteinver}
Y.~H. Cheng, T.~E. Wang, C.~J. Wu, A note on eigenvalue asymptotics for
  \textrm{Hill}'s equation., Appl. Math. Lett. 23~(9) (2010) 1013--1015.

\bibitem{Poschel}
J.~P{\"{o}}schel, E.~Trubowitz, Inverse Spectral Theory, Academic Press,
  Boston, 1987.

\bibitem{Melda.O}
O.~A. Veliev, M.~Duman, The spectral expansion for a nonself-adjoint
  \textrm{Hill} operator with a locally integrable potential, J. Math. Anal.
  Appl. 265 (2002) 76--90.

\bibitem{Veliev:Kiraç}
O.~A. Veliev, A.~A. \textrm{K\i ra\c{c}}, On the nonself-adjoint differential
  operators with the quasiperiodic boundary conditions, International
  Mathematical Forum 2~(35) (2007) 1703--1715.

\end{thebibliography}
\end{document}